\documentclass[12pt]{article}
\usepackage{amscd}
\usepackage{amsmath}
\usepackage{amsfonts}
\usepackage{amssymb}
\input{epsf}

\newcommand{\bbC}{{\mathbb C}}

\newcommand{\bbN}{\mathbb{N}}
 
\newcommand\Vect{\mbox{Vect}}
\newcommand\Hom{\mbox{Hom}}

\newcommand\Aut{\mbox{Aut}}
\newcommand\GL{\mathcal{GL}}

\newcommand\G{{\mathcal G}}
\newcommand\action{{\vartriangleright}}

\newtheorem{theorem}{Theorem}[section]  %reset counter in new section
\newtheorem{Lem}[theorem]{Lemma}
\newtheorem{Thm}[theorem]{Theorem}

\newtheorem{definition}[theorem]{Definition}
\newtheorem{Def}{Definition}[section]

\newtheorem{examples}[theorem]{Examples}
\newtheorem{Ex}[theorem]{Example}

\newtheorem{Conv}[theorem]{Convention}
\newtheorem{Rem}[theorem]{Remark}

\def\block{\hbox{${\vcenter{\vbox{\hrule height 0.4pt\hbox{\vrule width 0.4pt height 6pt \kern 5pt\vrule width 0.4pt}\hrule height 0.4pt }}}$}}
\def\qed{\hfill\block\vskip.2cm}

\newcommand\A{{\mathcal A}}

\newcommand\C{{\mathcal C}}
\newcommand\D{{\mathcal D}}
\newcommand\V{{\mathcal V}}
\newcommand\cG{{\mathcal G}}

\newcommand\cV{{\mathcal V}}

\newcommand\oG{\overline{\G}}

\renewcommand\th{\tilde{h}} % don't know what it replaces
\newcommand\tk{\tilde{k}}

\newcommand\To{\Rightarrow}

\newcommand\ol{\overline}

\newcommand\gpd{\,\lower3pt\hbox{$\longrightarrow$}\hskip-.267in\raise3pt
             \hbox{$\longrightarrow$}\,}

\newcommand\bicat{\mbox{bicat}} % bicat

\begin{document}

\title { Categorical representations  \\of categorical groups  }

\author{ John W. Barrett\\ 
School of Mathematical Sciences\\ University
Park\\ Nottingham NG7 2RD, UK \and 
Marco Mackaay\\
Departamento de Matem\'{a}tica\\ Universidade do Algarve\\
8005-139 Faro, Portugal}

\date{17th September 2004}  %overrides today's date
\maketitle 

\begin{abstract}  
The representation theory for categorical groups is constructed. Each categorical group determines a monoidal bicategory of representations. Typically, these categories contain representations which are indecomposable but not irreducible. A simple example is computed in explicit detail. 
\end{abstract}

\section{Introduction}
In three-dimensional topology there is a very successful interaction between category theory, topology, algebra and mathematical physics which is reasonably well understood. Namely, monoidal categories play a central role in the construction of invariants of three-manifolds (and knots, links and graphs in three-manifolds), which can be understood using quantum groups and, from a physics perspective, the Chern-Simons functional integral. The monoidal categories determined by the quantum groups are all generalisations of the idea that the representations of a group form a monoidal category.

The corresponding situation for four-manifold topology is less coherently understood and one has the feeling that the current state of knowledge is very far from complete. The complexity of the algebra increases dramatically in increasing dimension (though it might eventually stabilise). Formalisms exist for the application of categorical algebra to four-dimensional topology, for example using Hopf categories \cite{CF}, categorical groups \cite{Y-HT} or monoidal 2-categories \cite{CS,BL,M-S}. Since braided monoidal categories are a special type of monoidal 2-category (ones with only one object), then there are examples of the latter construction given by the representation theory of quasi-triangular Hopf algebras. This leads to the construction of the four-manifolds invariants by Crane, Yetter, Broda and Roberts which give information on the homotopy type of the four-manifold \cite{CKY,R,R-EX}. At present it seems that categorical invariants which delve further into the smooth or combinatorial structure of four-manifolds will require different types of examples of monoidal 2-categories.

In this paper we determine a new set of examples of monoidal 2-categories. We show that the categorical representations of a categorical group form a monoidal 2-category, by direct analogy with the way in which the representations of a group form a monoidal category. The categorical definitions are given in section 2 and 3. In section 4, an abstract definition of categorical representations, and their morphisms, is given and then unpacked. An extended example is calculated in section 4 with explicit matrices, illustrating many of the complexities of more general examples. In section 5 we give a fairly complete characterisation of the one-dimensional categorical representations, outlining a number of examples. Finally in section 6 we make a number of remarks about the structure of the N-dimensional categorical representations, in particular the phenomenon of representations which are indecomposable but not irreducible. These remarks generalise some of the features of the example in section 4. 

One particular example of a monoidal category leads to a state-sum model for quantum gravity in 
three-dimensional space-time \cite{TV,B-O}. The motivation for this work grew out of wondering if there is a corresponding model in the more realistic four-dimensional space-time. The first attempts at doing this \cite{BC,DFKR} used a braided monoidal category and suffer from several problems, one of which is that there is no `data' on the edges of a triangulation of the manifold, which is where one might expect to find the combinatorial version of the metric tensor \cite{REG}.  Thus we arrived at the idea of constructing the monoidal 2-category of representations for the example of the categorical Lie group determined by the Lorentz group and its action on the translation group of Minkowski space, generalising the construction of \cite{M-FG}. An early draft of this paper is the reference cited by Crane and Yetter \cite{CY-2G,CY-MC,Y-MC,CSH} who developed the particular example, and the machinery of measurable categories to handle the Lie aspect, much further.

\section{Categorical groups}
\label{catgroups}
\begin{definition} A {\em categorical group} is by definition a group-object 
in the category of groupoids.
\end{definition}
This means that a categorial group is a groupoid $\G$, with a set of objects $\G_0\subset\G$,
together with functors which implement the group product, $\circ\colon\G\times\G\to\G$, 
and the inverse $\mathstrut^{-1}\colon\G\to\G$, together with an identity object $1\in \G_0$. These satisfy 
the usual group laws:
 $$a\circ (b\circ c)= (a\circ b)\circ c$$
 $$a\circ 1=1\circ a=a$$
 $$a\circ a^{-1}=a^{-1}\circ a=1$$
for all $a$, $b$, $c$ in $\G$. In particular, $\G$ is a strict monoidal category.

\begin{definition}A {\em functorial homomorphism} between two categorical 
groups is a strict monoidal functor.  
\end{definition} 

Categorical groups are equivalent to crossed modules of groups. 
This equivalence, and the basic properties of categorical groups, 
are explained in \cite{BS}. Here we give a brief outline.

In a categorical group $\G$ with hom-sets $\G(X,Y)$, the categorical composition $f\cdot g$ and the group product $\circ$ are related by the interchange law
$$(f\circ g)\cdot(h\circ k)=(f\cdot h)\circ (g\cdot k).$$
For all categories in this paper the diagrammatic order of composition is used. This means that $f\cdot g$ is defined when the target of $f$ is equal to the source of $g$. 

In fact the composition is determined by the product. If $f\in\G(X,Y)$ and $g\in\G(Y,Z)$, then
$$f\cdot g=f\circ 1_{Y^{-1}}\circ g=g\circ 1_{Y^{-1}}\circ f.$$ In particular, 
$\G(1,1)$ is an abelian group, the composition and product coinciding.

\begin{definition}A {\em crossed module} is a homomorphism of groups
$$\partial \colon E\to G$$
together with an action $\action$ of $G$ on $E$ by automorphisms, such that
$$ \partial \bigl( X\action e\bigr) = X\bigl(\partial e\bigr)X^{-1}$$
$$ \bigl(\partial e\bigr)\action e' = e e' e^{-1}.$$
We call $E$ the {\em principal group} and $G$ the {\em base group}. 
\end{definition}

There is a natural notion of a mapping between crossed modules. 

\begin{definition}
\label{homcross}
A {\em homomorphism of crossed modules} 
$$(E,G,\partial,\action)\to(E',G',\partial',\action')$$
is given by two vertical homomorphisms 
$$\begin{CD}
E @>\partial>> G\\
@V{F_p}VV @VV{F_b}V \\
E' @>\partial'>> G'
\end{CD}$$
which commute and satisfy $F_p(X\action e)=F_b(X)\action'F_p(e)$. We call 
the latter condition the {\em action condition}. 
\end{definition}

The equivalence with categorical groups is as follows.

\begin{Thm}(Verdier)\label{Verdier} The category of categorical groups and 
functorial homomorphisms and the category of crossed modules of groups and 
homomorphisms between them are equivalent. 
\end{Thm}

The proof is sketched. Given a categorical group, 
a crossed module is defined by taking  the base group $G(\cG)$ 
to be the objects, 
$$G(\cG)=(\G_0,\circ),$$ 
and the principal group $E(\cG)$, to be the subset of morphisms which are 
morphisms from the object $1$ (to any object), 
$$E(\cG)=\bigcup_X \G(1,X),$$ 
again with the product operation $ab=a\circ b$. 
The homomorphism $E\to G$ of the crossed module is $e\in \G(1,X)\mapsto X$ and the action is $Y\action e=1_Y\circ e\circ 1_{Y^{-1}}$.

Conversely, given a crossed module $(E,G,\partial,\action)$, a categorical group is constructed in a 
canonical way by taking $\G_0=G$, $\G(X,Y)=\partial^{-1}\bigl(YX^{-1}\bigr)$. If $f\in\G(X,Y)$ and 
$g\in\G(Z,T)$, then the tensor product is defined as 
$$f\circ g = f(X\action g)\in\G(XZ,YT)$$
and the composition, for $Y=Z$,
$$ f\cdot g = gf\in\G(X,T).$$
\qed

It is worth noting that the 
definition of a crossed module and of a categorical group makes sense 
when the groups and groupoids are Lie groups and the equivalence between 
Lie categorical groups and Lie crossed modules holds in the same way.

\begin{examples}

\begin{enumerate}
\item \label{groupclosure} Let $K$ be a group, then we define $\overline{K}$, the {\em closure} of $K$, to be the groupoid with one object, 
$\bullet$, and hom-space $\overline{K}(\bullet,\bullet)=K$. If $K$ is abelian, then the group operation also defines a monoidal structure, so that $\overline{K}$ becomes a categorical group.
The corresponding 
crossed module has principal group $E(\overline{K})=K$ and trivial 
base group.   
\item A categorical group is {\em transitive} if there is a morphism 
between any pair of objects. The corresponding crossed module 
$\partial\colon E\to G$ is surjective. This implies that $E$ is 
a central extension of $G$, and the action of $G$ on $E$ is determined 
by conjugation in $E$.
\item A categorical group is {\em intransitive} if there are no 
morphisms between distinct objects (generalising (\ref{groupclosure})). In this case the crossed module 
$\partial\colon E\to G$ has $\partial =1$ and is determined entirely by 
the action of $G$ on the abelian group $E$. An example is the {\em wreath product}, 
where $G$ is the symmetric group $S_n$ and $E=(\bbC^*)^n$, 
with the action of $S_n$ being the permutation of the factors in $E$. 
%This example is the key example for categorical representations as we show in Lem.~\ref{}.

\item A categorical group is {\em free} if there is at most one 
morphism between any pair of objects. The crossed module  
$\partial\colon E\to G$
is injective and $E$ is a normal subgroup of $G$. Again the action of 
$G$ on $E$ is determined by conjugation.
\item The {\em transformation categorical group} of a category 
$\mathcal C$ is defined by 
$$\G_0=\{\mbox{functorial isomorphisms}\}\quad\mbox{and}\quad
\G=\{\mbox{natural isomorphisms}\}.$$ 
Here we consider, as usual, $\G_0$ to be a subset 
of $\G$ by identifying a functorial isomorphism with the identity natural isomorphism 
of it. For example, if $\mathcal C=\overline{K}$, 
then the crossed module corresponding to the transformation categorical group 
is $K\to \Aut K$, where $\partial$ maps a group element to the corresponding 
inner automorphism.
\end{enumerate}
\end{examples}

\section{Bicategories}

In this section we recall the definitions of $2$-dimensional category theory
\cite{G}. First we define 2-categories, sometimes called strict 2-categories, and then indicate the changes required to give the weaker notion of bicategories. Finally we discuss monoidal structures on bicategories and the example of 2-Vect.

\begin{Def} A {\em $2$-category}, $\C$, is given by:
\begin{enumerate}
\item A set of {\rm objects}, $\C_0$. 
\item A small category, $\C(X,Y)$, for each pair $X,Y\in\C_0$. 
The set of objects 
in $\C(X,Y)$ we denote by $\C_1(X,Y)$. The elements of $\C_1(X,Y)$ are called 
{\rm $1$-morphisms}. For each pair $f,g\in\C_1(X,Y)$, we denote the set of 
morphisms in 
$\C(X,Y)$ from $f$ to $g$ by $\C_2(f,g)$. The elements of $\C_2(f,g)$ are 
called {\rm $2$-morphisms}. The composition in $\C(X,Y)$ is called 
the {\rm vertical composition} and denoted by a small dot, e.g. 
$\mu\cdot\nu$ (or sometimes by simple concatenation without dot). 
\item A functor
$$\circ\colon \C(X,Y)\times\C(Y,Z)\to\C(X,Z),$$ 
for each triple $X,Y,Z\in\C_0$. Together these are 
called the {\rm horizontal composition}. The horizontal composition is 
required to be associative and unital. The last condition means that there 
is a $1$-morphism $1_X\in\C_1(X,X)$, 
for each $X\in\C_0$, which is a right and left unit for horizontal composition. 
\end{enumerate} 
\end{Def}

\begin{Ex}
\label{closure}
Let $\G$ be a categorical group. Then we define the {\rm closure 
of $\G$}, which we denote by $\oG$, to be the $2$-category with one 
object, denoted $\bullet$, such that $\oG(\bullet,\bullet)=\G$. The horizontal 
composition is defined by the monoidal structure in $\G$. Note that there is a slightly 
confusing mixture of subscripts now, because $\oG_0=\{\bullet\}$, 
whereas $\oG_1=\oG_1(\bullet,\bullet)=\G_0$. Unfortunately this renumbering seems unavoidable in 
this subject and we hope that the context will always avoid confusion in this paper.  
\end{Ex}
 
Next we recall the definition of $2$-functors, 
natural $2$-transformations and modifications.

\begin{Def} Given two $2$-categories, $\C$ and $\D$, a {\em $2$-functor} 
between 
them, $F\colon\C\to\D$, consists of:
\begin{enumerate}
\item A function $F_0\colon\C_0\to\D_0$.
\item A functor $F_1(X,Y)\colon\C(X,Y)\to\D(F(X),F(Y))$, 
for each pair $X,Y\in\C_0$. 
These functors are required to be compatible with the horizontal composition 
and the unit $1$-morphisms. By abuse of notation we sometimes denote both $F_0$ and 
$F_1(X,Y)$ simply by $F$.
\end{enumerate}
\end{Def}

\begin{Def}\label{nattwo} Given two $2$-functors between two $2$-categories, $F,G\colon\C\to\D$, a 
{\em natural $2$-transformation} between 
$F$ and $G$, denoted $h\colon F\To G$, consists of:
\begin{enumerate}
\item A $1$-morphism $h(X)\in\D_1(F(X),G(X))$, for each $X\in C_0$.
\item For each pair $X,Y\in\C_0$, a natural isomorphism 
$\th\colon F_1(X,Y)\circ h(Y)\to h(X)\circ G_1(X,Y)$, where we consider 
$F_1(X,Y)\circ h(Y)$ and $h(X)\circ G_1(X,Y)$ both as functors from 
$\C(X,Y)$ to $\D(F(X),G(Y))$. We 
require two coherence conditions to be satisfied: 
\begin{enumerate}
\item  $\th(f\circ g)=(1_{F(f)}\circ \th(g))\cdot(\th(f)\circ 1_{G(g)})$.
\item $\th(1_X)=1_{h(X)}$.
\end{enumerate}
\end{enumerate}
\end{Def}
 
This notion of natural $2$-transformation is almost the same as the notion of {\em quasi-natural transformation} in \cite{G}, which differs in that the 2-cells are not required to be isomorphisms and have the arrows reversed.

\begin{Def}\label{horcomnat2} Let $F,G,H\colon\C\to\D$ be three $2$-functors and let 
$h\colon F\To G$ and $k\colon G\To H$ be two natural $2$-transformations. 
The {\em horizontal composite} of $h$ and $k$, denoted $h\circ k$, 
is defined by 
\begin{enumerate}
\item For each $X\in \C_0$, $(h\circ k)(X)=h(X)\circ k(X)$.
\item For each pair $X,Y\in\C_0$, $\widetilde{h\circ k}=(\th\circ 1_{k(Y)})\cdot(1_{h(X)}\circ\tk)$. 
\end{enumerate}
\end{Def}

\begin{Def} Let $F,G\colon\C\to\D$ be two $2$-functors and 
$h,k\colon F\To G$ be two natural $2$-transformations. A {\em modification} 
$\phi\colon h\Rrightarrow k$ is given by a $2$-morphism $\phi(X)\in\D_2(h(X),k(X))$, 
for each $X\in\C_0$. These $2$-morphisms are required to satisfy
$$\th(f)\cdot(\phi(X)\circ 1_{G(f)})=(1_{F(f)}\circ \phi(Y))\cdot \tk(f),$$
for any $X,Y\in\C_0$ and $f\in \C_1(X,Y)$.   
\end{Def}

\begin{Def} The {\em horizontal} and {\em vertical compositions} of 
modifications are directly induced by the corresponding compositions in 
$\D$.  
\end{Def}
 
\subsection{Weakenings} The notion of a 2-category can be weakened to the notion of a bicategory, in which the associativity of horizontal composition is not given by an equation between 1-morphisms, but by 2-isomorphisms which are the components of a natural isomorphism called the {\em associator}, 
$$\alpha_{X,Y,Z}\colon f\circ (g\circ h)\To (f\circ g)\circ h,$$
which satisfies a certain coherence law. Similarly, the unital nature of the horizontal composition can be weakened by introducing natural isomorphisms \cite{BEN}. However, in all the examples in this paper the unital isomorphisms are all identities. We call this sort of bicategory a {\em strictly 
unital bicategory} 

The notions of 2-functor, natural $2$-transformation and modification have suitable generalisations to the case of bicategories, see \cite{ROU}, for example. One new phenomenon which occurs is that, for the 2-functors, the horizontal composition is no longer preserved exactly, but only up to a family of natural isomorphisms, defined as follows. If $F$ is a 2-functor and $f$ and $g$ composable 1-mophisms in its domain, then the 2-isomorphisms are
\begin{equation}
\label{weakfunc}
{\tilde F}_{fg}\colon F(f)\circ F(g) \to F(f\circ g).  
\end{equation}
In general a 2-functor between bicategories involves weakening the condition $F(1_f)=1_{F(f)}$ to an isomorphism. This notion of 2-functor is called a {\em homomorphism} in \cite{BEN, STR}. 

\begin{Def} A {\em strictly unitary homomorphism} is a 2-functor between bicategories 
for which the isomorphisms $F(1_f)\cong 1_{F(f)}$ are all identity 2-morphisms. If additionally, for all $f,g$, $\tilde F_{fg}$ is an identity 2-morphism we call $F$ a 
{\em strict homomorphism}
\end{Def}

The natural 2-transformations have a straightforward generalisation. They are called natural pseudo-transformations in \cite{ROU}. 

Given two bicategories ${\mathcal C}$ and ${\mathcal D}$, then the following result gives a construction of a new bicategory \cite{STR,GPS}.

\begin{Thm} \label{bicat} The 2-functors from $\mathcal C$ to $\mathcal D$, 
together with their natural $2$-transformations and their modifications form a 
bicategory $\bicat({\mathcal C},{\mathcal D})$.
\end{Thm}

\subsection{Monoidal structures} 
The notion of a monoidal stucture on a 2-category is straightforward to define.

 \begin{Def} A {\em monoidal $2$-category} is a $2$-category, $\C$, together
with a $2$-functor $\boxtimes\colon\C\times\C\to\C$, the {\em monoidal product},
which is associative and unital. The latter means that there is a
{\rm unit object}, $I$, which is a left and right unit for the
monoidal product.
\end{Def}
By $\C\times\C$ we mean the cartesian product $2$-category. The requirement that
$\boxtimes$ is a 2-functor means that the interchange law is satisfied as an
identity. This is too strict for the purposes of this paper, as we need to use bicategories rather than 2-categories. This gives the notion of a {\em monoidal bicategory}, for which the definition is the same as the definition for a  monoidal $2$-category, but with the appropriate notion of $2$-functor for bicategories and natural 2-transformations which give the associative and unital conditions. This means the $2$-functor $\boxtimes\colon\C\times\C\to\C$ carries natural 2-isomorphisms $\widetilde\boxtimes$ called the {\em tensorator}. The other coherers of the monoidal structure in the definition of a monoidal bicategory will always be trivial in this paper. This entire structure, together with the axioms it obeys, is a special case of the definition of a tricategory given by \cite{GPS}, in which the tricategory has only one object.

\subsection{2-Vector spaces}
\label{2vect}
In this section we recall the definition, due to Kapranov and
Voevodsky~\cite{KV}, of the monoidal bicategory of
2-vector spaces in the completely coordinatized version.  

\begin{Def} We define the {\rm  monoidal bicategory $2\Vect$} as
follows:
\begin{enumerate}
\item $2\Vect_0=\bbN$, the set of natural numbers including zero.
\item For any $N,M\in 2\Vect_0$, we define the category $2\Vect(N,M)$ as 
follows:
\begin{description}
\item{i)} The set of objects of $2\Vect_1(N,M)$ consists of all 
$N\times M$ matrices with coefficients in $\bbN$.
\item{ii)} For any $a,b\in2\Vect_1(N,M)$, 
the set $2\Vect_2(a,b)$ consists of all $N\times M$ matrices whose 
coefficients are complex matrices such that the $i,j$-coefficient has 
dimension $a^i_j\times b^i_j$. 
For any $a\in 2\Vect_1(N,M)$, we define $1_a$ to be the 
$N\times M$ matrix such that $(1_a)^i_j$ is the identity matrix of
dimension $a^i_j$ 
The {\em vertical composition} of two composable $2$-morphisms 
is defined by componentwise matrix multiplication  
$$(\alpha\cdot\beta)^i_j=\alpha^i_j\beta^i_j.$$
\end{description}

\item The {\em unit $1$-morphism} on an object $N$, denoted $1_N$, 
is given by the $N\times N$ identity matrix. 
The {\em (horizontal) composition} of 
two $1$-morphisms is defined by matrix multiplication.

The {\em horizontal composition} of two 
composable 
$2$-morphisms is defined by 
$$(\alpha\circ\beta)^i_j=\oplus_k\alpha^i_k\otimes \beta^k_j.$$ 

This composition is not strictly associative, and so there are associativity isomorphisms which carry out the corresponding permutations of bases.

\item The {\em monoidal product} of two objects
is defined by 
$$N\boxtimes M=NM.$$ 
For two 
$1$-morphisms we define 
$$(a\boxtimes b)^{ij}_{kl}=a^i_k b^j_l.$$ 
Finally, for two $2$-morphisms we define  
$$(\alpha\boxtimes\beta)^{ij}_{kl}=\alpha^i_k\otimes\beta^j_l.$$ 
The {\em unit object} is $1$.

This monoidal product has a tensorator $\widetilde\boxtimes$ which is again given by the corresponding permutations of basis elements. 

\item The remaining coherers for a monoidal bicategory are trivial.

\end{enumerate}
\end{Def}

So far we have defined the structure of the monoidal $2$-category. There
are also linear structures in $2\Vect$, which we define below. Although we
have not defined linear structures in general, we hope that the definitions
below are clear.  

\begin{Def}
\label{boxplus}
 There are three levels of linear structure in $2\Vect$, which we call
the {\em monoidal sum}, the {\rm direct sum} and the {\rm sum} respectively.
\begin{enumerate}
\item
The {\em monoidal sum} defines a monoidal structure on $2\Vect$.
For two objects it is defined as
$$N\boxplus M=N+M,$$
for two $1$-morphisms as
$$a\boxplus b=a\oplus b$$
and for two $2$-morphisms as
$$\alpha\boxplus \beta=\alpha\oplus\beta.$$
The {\em zero object} is $0$. Unlike $\widetilde\boxtimes$, which is non-trivial, $\tilde\boxplus$ 
is just the identity.
\item
The {\em direct sum} defines a monoidal structure on each $2\Vect(N,M)$. 
The direct sum of two $1$-morphisms $a,b$ in $2\Vect(N,M)$ is defined by
$$(a\oplus b)^i_j=a^i_j + b^i_j.$$ 
The direct sum of two $2$-morphisms in $2\Vect(N,M)$ is defined by
$$(\alpha\oplus\beta)^i_j=\alpha^i_j\oplus \beta^i_j.$$
The {\em zero $1$-morphism}, denoted $(0)$, is given by the
$N\times M$-dimensional zero matrix.
\item
The {\em sum} defines the linear structure in each $2\Vect(N,M)$. For two 
$2$-morphisms we define
$$(\alpha+\beta)^i_j=\alpha^i_j+\beta^i_j.$$
The {\rm zero $2$-morphism}, denoted $((0))$, is given by the matrix all
of whose entries are zero matrices of the right size. 
\end{enumerate}
\end{Def}

\begin{Def} Let $N\in \bbN$. We define the {\em general linear categorical 
group}, $\GL(N)\subset 2\Vect(N,N)$, to be the categorical group consisting
of all invertible $1$- and $2$-morphisms in $2\Vect(N,N)$.
\end{Def}
This definition makes sense because the horizontal composition of invertible 1- and 2-morphisms in $2\Vect$ is in fact strictly associative, as proved in the following lemma.

\begin{Lem} 
\label{wreath}
The monoidal category $\GL(N)$ is a categorical group, and the associated crossed module of groups is given by 
$$
\begin{CD}
(\bbC^*)^N@>{1}>> S_N,
\end{CD}
$$
where $S_N$ is the symmetric group on $N$ letters and $1$ the trivial 
group homomorphism which maps everything to $1\in S_N$. The action of $S_N$
on $(\bbC^*)^N$ is given by the permutations of the coordinates.  
\end{Lem}
\noindent{\bf Proof} The only invertible $1$-morphisms in 
$2\Vect(N,N)$ are the permutation matrices. The associator for a triple of permutation matrices is trivial, hence the horizontal composition in $\GL(N)$ is strictly associative and it forms a categorical group.

Between two permutation matrices there can only be an invertible $2$-morphism if they are equal. The reason is that, if
the corresponding entries in two permutation matrices are different, then one 
of them has to be zero. 

This shows that
$2\Vect_2(1_N,1_N)$ has its only non-trivial entries on the diagonal
and these are $1$-dimensional invertible complex
matrices. Hence $2\Vect_2(1_N,1_N)$ is isomorphic to $(\bbC^*)^N$. Theorem~\ref{Verdier} shows that 
the action of a permutation matrix $P\in2\Vect_1(N,N)$ on a diagonal matrix
$a\in2\Vect_2(1_N,1_N)$ is given by 
$$PaP^{-1}.$$
Here we mean the ordinary matrix multiplication, which is what 
$1_P\otimes a\otimes 1_{P^{-1}}$ amounts to in this case.
Thus the element in $S_N$, which corresponds to $P$, acts on the element in 
$(\bbC^*)^N$, which corresponds to $a$, by permutation of
its coordinates.\qed

\section{Categorical representations}\label{seccatrep}

In this section we give the abstract definitions of categorical 
representations, the functors between them, 
the natural transformations between such functors and the monoidal product. These definitions are analogues of Neuchl's~\cite{N} definitions for Hopf categories. We 
also work out a concrete example. 

 Let $G$ be a group and $\ol{G}$ its closure (defined in Example \ref{groupclosure}). 
One can check that a representation of $G$ is 
precisely a functor $\ol{G}\to\Vect$ and an intertwiner precisely a natural 
transformation between two such functors. This observation motivates 
the following definition for categorical groups.   
 
The starting point is the selection of a monoidal bicategory as the category in 
which the categorical group is represented. In this paper we discuss only 2Vect, 
although other monoidal bicategories could be used instead (see the remarks in the 
final section).

Let $\G$ be an arbitrary categorical group and let $\oG$ be its closure, defined in Ex.~\ref{closure}. We first give a conceptual definition. 

\begin{Def} a) A {\em categorical representation} of $\G$ is a strictly unitary homomorphism 
$(R,{\tilde R})\colon\oG\to 2\Vect$. We call the non-negative integer $R(\bullet)\in 2\Vect_0$ the 
{\em dimension} of the categorical representation.  

b) A {\em $1$-intertwiner} is a natural $2$-transformation between two 
categorical representations. 

c) A {\em $2$-intertwiner} is a modification between two $1$-intertwiners 
with the same source and target.

\end{Def}

\noindent Note that Neuchl~\cite{N} uses the terms $\G$-functors and 
$\G$-transformations instead of $1$- and $2$-intertwiners.  

Theorem~\ref{bicat} shows that the categorical representations of $\G$, 
together with the $1$- and $2$-intertwiners, form a bicategory, $\bicat(\oG,2\Vect)$. 
To give this the stucture of a monoidal bicategory, we first promote $\oG$ and $2\Vect$ to tricategories and use some general results about those. We consider $\oG$ as a strict tricategory, denoted $\underline{\oG}$ by adding only identity 3-morphisms to the existing strict bicategory. For $2\Vect$ we take the closure $\overline{2\Vect}$ as the tricategory. By a general result of Gordon, Power and Street~\cite{GPS} about tricategories we know that $\mbox{tricat}(\underline{\oG},\overline{2\Vect})$ forms a 
tricategory. The objects of this tricategory are trihomomorphisms (functors between tricategories). The constant trihomomorphism is the one that sends every 3-morphism to $1_{1_{1_\bullet}}$.

\begin{Lem} $\overline{\bicat(\oG,2\Vect)}$ is the subtricategory of $\mbox{tricat}(\underline{\oG},\overline{2\Vect})$ determined by the 1-, 2- and 
3-morphisms on the unique constant trihomomorphism.
\end{Lem}

\noindent{\bf Proof} Just check the diagrams in \cite{GPS}.\qed
\noindent Unfortunately Gordon, Power and Street~\cite{GPS} do not give explicit definitions of the composition rules for the various morphisms in tricategories. Therefore we spell out the tensor product in $\bicat(\oG,2\Vect)$ below, which corresponds to the 
horizontal composition in $\mbox{tricat}(\underline{\oG},\overline{2\Vect})$. The first 
definition can also be found in \cite{ROU}.

\begin{Lem}
\label{boxtimes}
Let $(R,{\tilde R})$ and $(T,{\tilde T})$ be two categorical representations. Then 
the monoidal product $(R\boxtimes T,\widetilde{R\boxtimes T})$ of $(R,{\tilde R})$ 
and $(T,{\tilde T})$ is given by  
$$R\boxtimes T(X)=R(X)\boxtimes T(X)$$
and the following diagram  
$$\begin{CD}
R(X)R(Y)\boxtimes T(X)T(Y) @>{{\tilde {R}_{X,Y}}\boxtimes{\tilde {T}_{X,Y}}}>> R(XY)\boxtimes T(XY)=R\boxtimes T(XY)\\
@VV{\widetilde{\boxtimes}_{(R(X),T(X)),(R(Y),T(Y))}}V  @ AA{\widetilde{R\boxtimes T}_{X,Y}}A\\
(R(X)\boxtimes T(X))(R(Y)\boxtimes T(Y))@= (R\boxtimes T(X))(R\boxtimes T(Y))\\
\end{CD}$$
\end{Lem}

\begin{Lem}
Let $(h,\th)\colon R_1\to R_2$ and $(k,\tk)\colon T_1\to T_2$ be $1$-intertwiners. 
The monoidal product,  
$(h\boxtimes k,\widetilde{h\boxtimes k})$, is given by 
$$ h\boxtimes k (X) = h(X)\boxtimes k(X)$$
using, on the right, the monoidal product in 2Vect, and by the following diagram
$$\begin{CD}
R_1\boxtimes T_1(f)\circ h\boxtimes k(Y)@>{\widetilde{h\boxtimes  k}}>>(h\boxtimes k(X))\circ(R_2\boxtimes T_2(f))\\
@| @|\\
(R_1(f)\boxtimes T_1(f))\circ (h(Y)\boxtimes k(Y)) &{}& (h(X)\boxtimes k(X))\circ(R_2(f)\boxtimes T_2(f))\\
@V{\widetilde\boxtimes}VV @V{\widetilde\boxtimes}VV\\
(R_1(f)\circ h(Y))\boxtimes (T_1(f)\circ k(Y))@>\tilde{h}\boxtimes\tilde{k}>> 
(h(X)\circ R_2(f))\boxtimes (k(X)\circ T_2(f))\\ 
\end{CD}$$ 
\end{Lem}

\begin{Lem} Let $\alpha$ and $\beta$ be two $2$-intertwiners. The monoidal product is
given by 
$$\alpha\boxtimes\beta,$$
using the monoidal product in $2\Vect$.  
\end{Lem}

There is also a natural way of defining the monoidal sum of two 
categorical representations.

\begin{Def} Let $(R,\tilde{R})$ and $(T,\tilde{T})$ be two categorical 
representations. Their {\em monoidal sum} is defined by 
$$R\boxplus T(X)=R(X)\boxplus T(X),$$
and 
$$\tilde{R}\boxplus\tilde{T}(X,Y)=\tilde{R}(X,Y)\boxplus\tilde{T}(X,Y).$$
Note that the latter makes sense, because 
$$(R(X)\circ R(Y))\boxplus (T(X)\circ T(Y))=(R(X)\boxplus T(X))\circ (R(Y)\boxplus T(Y))$$
holds on the nose, i.e. the 2-isomorphism between both sides of the equation is the identity, as already 
noticed in definition~\ref{boxplus}. 
\end{Def}

\noindent Since $\widetilde{\boxplus}$ is trivial, the definition of their 
monoidal sum of 1- and 2-intertwiners is much simpler than that of their 
monoidal product.

\begin{Def} Let $(h,\th)\colon R_1\to R_2$ and $(k,\tk)\colon T_1\to T_2$ be 
$1$-intertwiners. We define their {\em monoidal sum},  
$(h\boxplus k,\widetilde{h\boxplus k})$, to be 
 
$$ h\boxplus k (X) = h(X)\boxplus k(X)$$
and  
$$\tilde{h}\boxplus\tilde{k},$$ 
using on the right the monoidal sum in 2Vect.
\end{Def}
\begin{Def}
Let $\alpha$ and $\beta$ be two $2$-intertwiners. We define 
their {\em monoidal sum} as
$$\alpha\boxplus\beta,$$
using the monoidal product in $2\Vect$. 
\end{Def}

\noindent We do not give a precise definition of the direct sum of 1- and 
2-intertwiners and the sum of 2-intertwiners, because we do not need them. 

\subsection{Strict categorical representations}
These definitions can be unpacked by applying the general definitions of bicategories 
to this particular situation and expressing the result in terms of categorical groups. 
The categorical representations can be formulated in terms of crossed modules. In 
this way a strictly unitary homomorphism leads to a weakened notion of morphism 
between crossed 
modules, involving a group 2-cocycle on the object group $G$ corresponding to 
$\tilde R$. The extra conditions on this cocycle appear to be very complicated. 
Therefore we restrict our attention to a subclass of categorical representations in the 
rest of this paper, which we call {\em strict categorical representations}.

\begin{Def} A {\em strict categorical representation} is a strict homomorphism 
$R\colon \oG\to 2\Vect$. This means that $\tilde R_{X,Y}$ are all identity 2-morphisms. 
\end{Def}      

 Restricting to strict categorical representations still gives a monoidal bicategory, as proved by the following lemma.

\begin{Lem} The monoidal product and the monoidal sum of two strict categorical representations yield 
strict categorical representations.  
\end{Lem}

\noindent{\bf Proof} The first claim follows from the fact that the top and the left-hand side of the diagram 
in lemma~\ref{boxtimes} are trivial. Note that 
$\widetilde{\boxtimes}_{(R(X),T(X)),(R(Y),T(Y))}$ is trivial, because all 
matrices in the subscript are permutation matrices. 

The second claim is obvious as well, because the monoidal sum of two identity 2-morphisms is an identity 
2-morphism. \qed

Let us now unpack the definition of a strict categorical representation and 1- and 
2-intertwiners.
 
\begin{itemize}
\item
A strict categorical representation of $\G$ amounts to a choice of a non-negative integer $N$ and a strict 
homomorphism $R\colon\oG\to\overline{\GL(N)}$, or, equivalently, a functorial homomomorphism between 
categorical groups $R\colon\G\to\GL(N)$.
It can also be described as a homomorphism between the corresponding crossed modules, 
\begin{equation}
\label{catrepcross}
\begin{CD}
E(\G) @>\partial>> G(\G)\\
@V{R_p}VV @VV{R_b}V \\
(\mathbb{C}^*)^N @>1>> S_N
\end{CD}
\end{equation}
according to theorem~\ref{Verdier} and lemma~\ref{wreath}.

\item
Let $R\colon\G\to\GL(N)$ and $T\colon\G\to\GL(M)$ 
be two strict categorical representations of $\G$. 
A $1$-intertwiner between them consists of a $1$-morphism 
$h_\bullet\in 2\Vect_1(N,M)$ together with a $2$-isomorphism $\th(X)\in 
2\Vect_2(R(X)\circ h_\bullet,h_\bullet\circ T(X))$, for each $X\in \G_0$.
 
The 1-morphism $h_\bullet$ can be thought of as a vector bundle (with fibres of varying dimension) over the finite set $N\times M$, the cartesian product of the $N$-element set with the $M$-element set. The data above determine an action of the group $\G_0$ on this vector bundle in the following way.

There is a right action of $X\in\G_0$ on the set $N$, $i\mapsto iR_b(X)$, given by the permutation matrix $R_b(X)$ (acting on the right) and similarly a right action of $\G_0$ on the set $M$ given by $T_b$. 

\begin{Lem}\label{1intunpack} The collection of linear maps
$$ \th(X)^i_j\colon {h_\bullet}^{iR_b(X)}_j\to {h_\bullet}^i_{jT_b(X)^{-1}} $$
determines a left action of $\G_0$ on the vector bundle $h_\bullet$.
\end{Lem}

\noindent{\bf Proof}
 The two coherence conditions in definition \ref{nattwo} are 
$$\th(X\circ Y)=\left(1_{R(X)}\circ\th(Y)\right)\cdot\left(\th(X)\circ 1_{T(Y)}\right)$$
$$\th(1)=1_{h_\bullet}.$$ Note that the associators in these expressions are trivial because $R$ and $T$ take values in the permutation matrices. Taking the $(i,j)$-th component of $\th$ gives
$$ \th(X\circ Y)^i_j=\th(Y)^{iR_b(X)}_j\cdot \th(X)^i_{jT_b(Y)^{-1}},$$
which is the condition for a left action.
This covers the left action of $\G_0$ on $N\times M$ given by $(iR_b(X),j)\mapsto (i,jT_b(X)^{-1})$.\qed

The remaining condition satisfied by the 1-intertwiner is the naturality condition. This is that the diagram
 
\begin{equation}
\label{Gfunctor1}
\begin{CD}
R(X)\circ h_\bullet@>\th(X)>> h_\bullet\circ T(X)\\
@V{R(f)\circ 1_{h_\bullet}}VV @VV{1_{h_\bullet}\circ T(f)}V\\
R(Y)\circ h_\bullet@>\th(Y)>> h_\bullet\circ T(Y)
\end{CD}
\end{equation}
commutes. 

Since $\th(X)$ is invertible, this implies the equation on 1-morphisms
\begin{equation}\label{Gfunctor3} R(X)\circ h_\bullet=h_\bullet\circ T(X).\end{equation}
Also, restricting naturality to the crossed module data ($X=1$), gives the equation
\begin{equation}
 \left(R(e)\circ 1_{h_\bullet}\right)\cdot \th(\partial e)=1_{h_\bullet}\circ T(e)\label{Gfunctor2}
\end{equation}
for $e\in E(\G)$.

It is possible to show that this last condition is actually sufficient to recover all of (\ref{Gfunctor1}).

\begin{Lem}\label{naturalitylemma}
For intertwiners of strict categorical representations,
the naturality condition (\ref{Gfunctor1}) follows from condition
(\ref{Gfunctor2})
 for all $e\colon 1\to Y$ and the coherence condition 2(a) of definition \ref{nattwo}.
\end{Lem}

\noindent{\bf Proof}
The naturality condition
follows from the following computation. Suppose $f\colon X\to X\circ Y$ is
an arbitrary morphism in $\G_1$. We can write $f=1_X\circ e$, where $e\colon 1\to
Y$. Then
\begin{align*}\left(R(f)\circ 1_{h_\bullet}\right)&\cdot\th(X\circ Y)\\
&=\left( 1_{R(X)}\circ R(e) \circ 1_{h_\bullet} \right)\cdot
   \left((1_{R(X)}\circ\th(Y)\right) 
 \cdot
   \left(\th(X)\circ 1_{T(Y)}\right)\\
  \intertext{(using 2(a) of definition \ref{nattwo})}
&= \left( 1_{R(X)}\circ
   \left(\left(R(e)\circ 1_{h_\bullet}\right)\cdot \th(Y)\right)\right)
   \cdot \left(\th(X)\circ 1_{T(Y)}\right)\\
&= \left(1_{R(X)}\circ 1_{h_\bullet}\circ T(e)\right)
   \cdot \left(\th(X)\circ 1_{T(Y)}\right)\\
\intertext{(using (\ref{Gfunctor2}))}
&=\th(X)\circ T(e)\\
&=\left(\th(X)\circ 1_1\right)\cdot
  \left(1_{h_\bullet}\circ 1_{T(X)}\circ T(e)\right)\\
&=\th(X)\cdot \left( 1_{h_\bullet}\circ T(f)\right).\end{align*}
In the computation, all associators which occur are the identity.\qed

\item
Let $h=(h_\bullet,\th),k=(k_\bullet,\tk)\colon R\to T$ be $1$-intertwiners. A 
$2$-intertwiner between them consists of a single $2$-morphism 
$\phi\in 2\Vect_2(h_\bullet,k_\bullet)$.

The condition that this satisfies is that for each $X\in\G_0$ the following diagram commutes:
\begin{equation}
\label{Gtransformation1}
\begin{CD}
R(X)\circ h_\bullet@>\th(X)>> h_\bullet\circ T(X)\\
@V{1_{R(X)}\circ\phi}VV  @VV{\phi\circ 1_{T(X)}}V\\
R(X)\circ k_\bullet @>\tk(X)>> k_\bullet\circ T(X)
\end{CD}
\end{equation}
\end{itemize}

\subsection{Example}

\begin{Conv}
From now on a categorical representation will always mean a {\em strict} categorical representation.  
\end{Conv}

In this section we work out the categorical representations of a concrete 
example of a finite intransitive categorical group. Recall that the crossed modules 
corresponding to intransitive categorical groups are determined by a group $G$ and an 
abelian group $E$ on which $G$ acts by automorphisms. The simplest example is 
$G=C_2$, the cyclic group with two elements $\pm1$, and 
$E=C_3=\{1,x,x^{-1}\}$, with the non-trivial action of $C_2$ on $C_3$, $-1\action x=x^{-1}$. We call 
this categorical
group $\G(2,3)$. Another way of looking at this example is by defining the total space 
of all the morphisms  
$$\G(2,3)=$$
$$\left\{\begin{pmatrix} 1& 0\\ 0& 1 \end{pmatrix},
\begin{pmatrix} 0& 1\\ 1& 0 \end{pmatrix},
\begin{pmatrix} x& 0\\ 0& x^{-1} \end{pmatrix},
\begin{pmatrix} x^{-1}& 0\\ 0& x \end{pmatrix},
\begin{pmatrix} 0& x\\ x^{-1}& 0 \end{pmatrix},
\begin{pmatrix} 0& x^{-1}\\ x& 0 \end{pmatrix}\right\}.$$
The source and target of each morphism is determined by the place of the non-zero 
coefficients of the corresponding matrix, e.g. the sources and targets of the first 
two matrices are equal to $1$ and $-1$ respectively. In this matrix notation the 
monoidal product, corresponding to the horizontal composition in $\overline{\G(2,3)}$ 
and denoted by $\circ$, is defined by matrix multiplication (giving the dihedral group $D_3$) 
and the composition, corresponding to the vertical composition in $\overline{\G(2,3)}$ 
and denoted by simple concatenation, by coefficientwise multiplication. 

The constructions will be carried out in this example. The main features of the general case are 
apparent in this example; some comments on the generalisations are given at the end of the section.

We classify all 1- and 2-dimensional categorical representations of $\G(2,3)$. 

\begin{enumerate}
\item $\V(1)$. This is the identity representation defined by $R_b(\pm 1)=1\in S_1$ and $R_p(x)=1\in\bbC^*$. It is 
the only 1-dimensional categorical representation, due to the following argument. Obviously $R_b(\pm 1)=1\in S_1$ has to hold. 
By (\ref{catrepcross}) we see that $R(x)=\xi(x)$, 
where $\xi$ is a complex group character on $C_3$. By the action condition in definition~\ref{homcross} we 
see that $\xi(x^{-1})=\xi(x)$, so $\xi$ has to be the 
trivial character and $R_p(x)=1$.
\item $\V(2)$. This is the trivial 2-dimensional categorical representation defined by 
$R_b(\pm 1)=1\in S_2$ and $R_p(x)=(1,1)\in(\bbC^*)^2$. As in the 1-dimensional case, the action condition forces $R_p$ to be trivial 
if $R_b$ is trivial. Below we show that $\V(2)$ is isomorphic to $\V(1)\boxplus\V(1)$. 
\item $\V(2)_{\xi}$. These are the non-trivial 2-dimensional categorical 
representations, where 
$R_b(\pm 1)=\pm 1\in S_2$ and $R_p$ is determined by one complex group character, $\xi$, on $C_3$. 

By (\ref{catrepcross}) we see that $R_p(x)=(\xi(x),\psi(x))$, where $\xi$ and $\psi$ are both complex 
group characters on $C_3$. The action condition now becomes $(\xi(x)^{-1},\psi(x)^{-1})=
(\psi(x),\xi(x))$, so we have $\xi=\psi^{-1}$. There are no further restriction on $\xi$. A nice 
way of picturing the strict homomorphisms $R$ of these representations is by using the 
matrix definition of $\G(2,3)$:
$$R\left(\begin{pmatrix}0&1\\1&0\end{pmatrix}\right)=\begin{pmatrix}0&1\\1&0
\end{pmatrix},\quad R\left(\begin{pmatrix}x&0\\0&x^{-1}\end{pmatrix}\right)=
\begin{pmatrix}\xi(x)&0\\0&\xi(x)^{-1}\end{pmatrix}.$$
The image of the other endomorphisms is obtained via horizontal composition.   
\end{enumerate}

Next we study all the 1-intertwiners between these categorical 
representations. In Lemma~\ref{1intunpack} we showed that these can be 
seen as vector bundles with a left action of $\G_0$. In particular, equation (\ref{Gfunctor3}) holds. The naturality condition simplifies a bit further in this example. Let $(h_\bullet,\tilde{h})$ be a 1-intertwiner between two categorical 
representations $R$ and $T$. By lemma \ref{naturalitylemma} the naturality condition 
reduces to the equation
\begin{equation}
\label{naturality}
 R(e) \circ 1_{h_\bullet} =1_{h_\bullet}\circ T(e),
\end{equation} 
for any $e\in E(\G(2,3))$. Note that this shows that the action on the vector bundle does 
not have to satisfy any additional conditions. 

First we study the 1-intertwiners between categorical representations of 
the same dimension.  
\begin{enumerate}
\item $\C(\V(1),\V(1))$: A $1$-intertwiner, 
in this case, is given 
by $h_\bullet=(n)$, i.e. a 1-dimensional matrix with a non-negative integer 
coefficient, and an $n$-dimensional representation of $C_2$ denoted by $\th$. The naturality condition for 1-endomorphisms does not impose any restrictions in this case, as one can easily check.   
\item $\C(\V(2),\V(2))$: A $1$-intertwiner is given by 
$$h_\bullet=\begin{pmatrix}n_1 & n_2\\ n_3 & n_4\end{pmatrix},$$
with $n_i\in\bbN$, for $i=1,2,3,4$, and $\tilde{h}$. The naturality condition (\ref{naturality}) does not 
impose any restrictions in this case, because $R$ is trivial. Just as in the one-dimensional case we see 
that $\tilde{h}$ defines a representation of $C_2$ on $\bbC^{n_i}$, for all $i=1,2,3,4$.   
\item $\C(\V(2)_{\xi},\V(2)_{\psi})$: In this case a 
$1$-intertwiner is given by 
$$
h_\bullet=\begin{pmatrix}n_1&n_2\\ n_3&n_4\end{pmatrix}
$$
and $\tilde{h}$. The naturality condition imposes restrictions on $h$. We need consider only equations (\ref{Gfunctor3}) and (\ref{naturality}). First,
let $X=\begin{pmatrix}0 &1\\1&0\\\end{pmatrix}$. This gives 

\begin{eqnarray}
\nonumber
\begin{pmatrix}
 1_{n_2}&  1_{n_1}\\  1_{n_4} &  1_{n_3}
\end{pmatrix}&=&
\begin{pmatrix}
1_{n_1}& 1_{n_2}\\ 
1_{n_3}& 1_{n_4}
\end{pmatrix}\circ 
\begin{pmatrix}
0&1\\ 
1&0
\end{pmatrix}\\
\label{1intcond}&=&
\begin{pmatrix}
0&1\\ 
1&0
\end{pmatrix}\circ
\begin{pmatrix}
1_{n_1}& 1_{n_2}\\ 
1_{n_3}& 1_{n_4}
\end{pmatrix}\\
\nonumber &=&
\begin{pmatrix}
 1_{n_3}& 1_{n_4}\\
 1_{n_1}& 1_{n_2}
\end{pmatrix}.
\end{eqnarray}
This shows that $n_1=n_4=n$ and $n_2=n_3=m$.

Now consider $e=\begin{pmatrix} x&0\\0& x^{-1}\end{pmatrix}$ in equation (\ref{naturality}). In the same way, this gives
$$ 
\begin{pmatrix}\xi(x)1_n & \xi(x)1_m\\ \xi(x)^{-1}1_m & \xi(x)^{-1}1_n \end{pmatrix}
=\begin{pmatrix} \psi(x)1_n & \psi(x)^{-1} 1_m\\ \psi(x)1_m & \psi(x)^{-1}1_n\end{pmatrix}.$$
Therefore there are three possible cases: a) $\psi=\xi\not\equiv 1$, b)
$\psi=\xi^{-1}\not\equiv 1$ and c) $\psi=\xi\equiv 1$. 

\begin{description}
\item{a)} Since the character group of $C_3$ 
is isomorphic to $C_3$ we see that $\xi\ne\xi^{-1}$. 
Equation (\ref{1intcond}) holds if and only if  
$n\in\mathbb{N}$ is arbitrary, and $m=0$. Thus we have 
$$h_\bullet=\begin{pmatrix}n&0\\0&n\end{pmatrix}.$$
The condition that $\th$ is an action on the vector bundle defined by $h$ 
implies that we have
$$\th(1)=\begin{pmatrix}1_n&0\\0&1_n\end{pmatrix}\quad\mbox{and}\quad
\th(-1)=\begin{pmatrix}0&A\\A^{-1}&0\end{pmatrix},$$
where $A\in\mbox{GL}(n,\mathbb{C})$ is arbitrary. 
\item{b)} Equation (\ref{1intcond}) holds if and only if $n=0$, and 
 $m\in\mathbb{N}$ is arbitrary. Thus we have 
$$h_\bullet=\begin{pmatrix}0&m\\m&0\end{pmatrix}.$$
Again we can determine $\th$ explicitly:
$$\th(1)=\begin{pmatrix}0&1_m\\1_m&0\end{pmatrix}\quad\mbox{and}\quad
\th(-1)=\begin{pmatrix}A&0\\0&A^{-1}\end{pmatrix},$$
where $A\in\mbox{GL}(m,\mathbb{C})$ is arbitrary. 
Taking $m=1$ shows that $\V(2)_{\xi}$ and $\V(2)_{\xi^{-1}}$ are isomorphic. 
\item{c)} Equation (\ref{1intcond}) holds with no restriction on $n$ and $m$. 
Thus we have
$$h_\bullet=\begin{pmatrix}n & m\\m & n\end{pmatrix}.$$
Obviously we can decompose this matrix as 
$$\begin{pmatrix}n & m\\m & n\end{pmatrix}=\begin{pmatrix}n & 0\\0 & n\end{pmatrix}\oplus\begin{pmatrix}0 & m\\m & 0\end{pmatrix}.$$
As in the previous two cases we see that 
$$\th(1)=\begin{pmatrix}1_n&1_m \\ 1_m&1_n\end{pmatrix}\quad\mbox{and}\quad \th(-1)=
\begin{pmatrix}A & B \\ B^{-1} & A^{-1}\end{pmatrix},$$
with arbitrary $A\in\mbox{GL}(n,\mathbb{C})$ and $B\in\mbox{GL}(m,\mathbb{C})$.    
\end{description}
\end{enumerate}
Finally let us describe the 1-intertwiners between categorical representations 
of different type. 
\begin{enumerate}
\item Let $(h_\bullet,\tilde{h})\colon \V(1)\to\V(2)$ be a 1-intertwiner. Then $h_\bullet$ has the form $h_\bullet=\begin{pmatrix}
n & m\end{pmatrix}$, 
where $n,m\in\mathbb{N}$. The naturality condition does not impose any restrictions because both 
categorical representations are trivial. As before, $\tilde{h}$ simply defines two representations of 
$C_2$ of dimensions $n$ and $m$ respectively. 

Analogously we see that any 1-intertwiner $(h_\bullet,\tilde{h})\colon\V(2)\to\V(1)$ is of the form 
$$h_\bullet=\begin{pmatrix}n\\m\end{pmatrix},$$
with $\tilde{h}$ defining two representations of $C_2$ again. 
Taking $n=1,m=0$ and $n=0,m=1$ respectively, and $\tilde{h}$ the trivial representation in both 
cases, shows that $\V(2)$ is isomorphic to $\V(1)\boxplus \V(1)$. 
\item Let $(h_\bullet,\tilde{h})\colon\V(1)\to\V(2)_{\xi}$ be a 1-intertwiner. Again $h_\bullet$ has the form 
$h_\bullet=\begin{pmatrix}n&m\end{pmatrix}$. The naturality condition now imposes the following 
restriction:
$$
\begin{pmatrix} 1_n & 1_m\end{pmatrix}\circ \begin{pmatrix}0&1\\ 1
&0\end{pmatrix}=\begin{pmatrix} 1_m & 1_n \end{pmatrix}=
\begin{pmatrix} 1_n & 1_m\end{pmatrix}.
$$
This holds if and only if $n=m$.
Taking $e=\begin{pmatrix} x&0\\0& x^{-1}\end{pmatrix}$, we see that $\xi$ has to be equal to $1$. 
Just as before $\tilde{h}$ defines an action on a vector 
bundle over $2$ with fibre $\mathbb{C}^n$.

Likewise non-zero 1-intertwiners $(h_\bullet,\tilde{h})\colon\V(2)_{\xi}\to 
\V(1)$ can be seen to exist if and only if $\xi$ is trivial and the intertwiners are the 
transposes of the previous ones. 
\item Let $(h_\bullet,\tilde{h})\colon \V(2)\to \V(2)_{\xi}$ be a 1-intertwiner. We already know that $\V(2)\cong 
\V(1)\boxplus \V(1)$. Therefore this case reduces to the direct sum of the previous case.
Again by transposition we get the classification of all 1-intertwiners between $\V(2)_{\xi}$ and $\V(2)$. 
\end{enumerate}

The 2-intertwiners are very easy to describe. Given two 1-intertwiners between 
two categorical representations, we know that we can interpret them as 
homogeneous vector bundles by the above results. A 2-intertwiner between them can then be interpreted 
as a bundle map between these 
homogeneous vector bundles which commutes with the actions of $C_2$. This 
interpretation follows immediately from the above and diagram (\ref{Gtransformation1}).    

Let us now have a look at the monoidal product of the above categorical representations.

\begin{enumerate}
\item Clearly we have $\V(1)\boxtimes \V\cong \V$, for any categorical representation $\V$. 
Because we also know that $\V(2)\cong \V(1)\boxplus\V(1)$, we see that 
$\V(2)\boxtimes \V\cong \V\boxplus \V$.
\item We now study $\V(2)_{\xi}\boxtimes\V(2)_{\psi}$. The easiest way to understand this tensor product 
is by looking at 
$$\begin{pmatrix}\xi(x)&0\\0&\xi(x)^{-1}\end{pmatrix}\boxtimes
\begin{pmatrix}\psi(x)&0\\0&\psi(x)^{-1}\end{pmatrix}=$$
$$\begin{pmatrix}\xi(x)\psi(x)&0&0&0\\0&\xi(x)\psi(x)^{-1}&0&0\\
0&0&\xi(x)^{-1}\psi(x)&0\\
0&0&0&\xi(x)^{-1}\psi(x)^{-1}\end{pmatrix}
$$
and
$$\begin{pmatrix}0&\xi(x)\\\xi(x)^{-1}&0\end{pmatrix}\boxtimes
\begin{pmatrix}0&\psi(x)\\\psi(x)^{-1}&0\end{pmatrix}=$$
$$\begin{pmatrix}0&0&0&\xi(x)\psi(x)\\
0&0&\xi(x)\psi(x)^{-1}&0\\
0&\xi(x)^{-1}\psi(x)&0&0\\
\xi(x)^{-1}\psi(x)^{-1}&0&0&0\end{pmatrix}.
$$
It is now easy to check that 
$$h_\bullet=\begin{pmatrix}1&0&0&0\\0&0&0&1\\0&1&0&0\\0&0&1&0\end{pmatrix}$$
and 
$$\tilde{h}=\begin{pmatrix}(1)&0&0&0\\0&0&0&(1)\\0&(1)&0&0\\0&0&(1)&0\end{pmatrix}$$
define an invertible 1-intertwiner  
$$\V(2)_{\xi}\boxtimes\V(2)_{\psi}\to \V(2)_{\xi\psi}\boxplus\V(2)_{\xi\psi^{-1}}.$$  
Note that $\xi\psi$ or $\xi\psi^{-1}$ is trivial, for any choice of $\xi$ and $\psi$.
\end{enumerate}

\section{Categorical characters}

In this section we study the 1-dimensional (strict) categorical representations 
of an arbitrary categorical group $\G$.

\begin{Thm} 
\label{catchar}
(a) A one-dimensional categorical representation, $R$, of $\G$ 
is completely determined by a group character, $\xi_{R}$, on 
$E=E(\G)$ which is invariant under the action of $G=G(\G)$. 

(b) A $1$-intertwiner $(h_\bullet,\tilde{h})$ between two one-dimensional categorical representations, 
$R$ and $T$, is either zero or given by a representation $\tilde{h}\colon G\to \mbox{GL}(h_\bullet)$, such that 
\begin{equation}
\label{Gfunccondition}
\tilde{h}(X)=\xi_R(e)^{-1}\xi_T(e),
\end{equation}
for $X=\partial e$ and any $e\in  E$. The right-hand side of 
(\ref{Gfunccondition}) should be read as a scalar matrix of the right size 
for the equation to make sense. 
In particular, there exists no non-zero $1$-intertwiner 
if the restrictions of $\xi_R$ and $\xi_T$ to $\ker\partial\subset E$ are different. 
If $R=T$, then $\tilde {h}$ has to be trivial on 
$\partial(E)$. The composition of two $1$-intertwiners corresponds to the 
tensor product of the two respective representations of $G$.

(c) Suppose we have two $1$-intertwiners between the same pair of categorical 
representations, then a $2$-intertwiner between them is given by an 
ordinary intertwiner between the corresponding representations. 
The horizontal composition of two $2$-intertwiners corresponds to the tensor product 
of the respective ordinary intertwiners, 
whereas the vertical composition of two 
$2$-intertwiners corresponds to the ordinary (matrix) product of the two respective intertwiners.   
\end{Thm}
\noindent{\bf Proof} 
(a) The character is $\xi_R(e)=R(e)$ restricted to $e\in E$. The first claim follows immediately from diagram 
(\ref{catrepcross}). 

(b) Now consider a $1$-intertwiner $(h_\bullet,\tilde{h})$ between $R$ and $T$. 
In this particular case 
$h_\bullet\in\Vect_1(1,1)\cong\bbN$, so we can identify $h$ with a natural number. 
Just as in the previous section we see that $\tilde{h}$ defines a representation of $G$ 
on $\bbC^{h_\bullet}$. Condition (\ref{Gfunccondition}) expresses the naturality condition for $h$ 
and can be read off from diagram (\ref{Gfunctor2}).

The composite of two $1$-intertwiners is given by the tensor product of the 
respective representations of $G$, say $\tilde{h}_1$ and $\tilde{h}_2$, 
because, in the case of one-dimensional 
categorical representations, we have (see definition \ref{horcomnat2}) 
$$\left((\tilde{h}_1)\circ 1\right)\cdot\left(1\circ(\tilde{h}_2)\right)=
(\tilde{h}_1\otimes \tilde{h}_2),$$
where the $1$'s denote the identity matrices of the right size.  

(c) This follows directly from condition (\ref{Gtransformation1}). 
\qed   
\vskip5pt
\noindent We now apply our results to the examples in Sect.~\ref{catgroups}. 

\begin{Ex} Let $\G$ be a categorical group and $(E,G,\partial,\action)$ the 
corresponding crossed module of groups. A one-dimensional categorical 
representation is determined by a $G$-invariant group character on $E$, 
say $\xi$, and is denoted by $\cV_{\xi}$. We determine the Hom-categories 
between an arbitrary pair $\cV_{\xi}$ and $\cV_{\psi}$. Let $(h_\bullet,\tilde{h})$ be 
an arbitrary $1$-intertwiner between $\cV_{\xi}$ and $\cV_{\psi}$.   
\begin{enumerate}
\item Suppose $\G$ is transitive. Then theorem~\ref{catchar} says that 
  $\tilde{h}$ is a scalar 
representation of $G$, which is completely determined by $\psi\xi^{-1}$. 
Therefore, for any two characters $\xi$ and $\psi$ which coincide on the kernel of $\partial$, there is exactly one 1-intertwiner between 
$\V_{\xi}$ and $\V_{\psi}$. The $2$-intertwiners are just the 
ordinary intertwiners between these scalar representations of $G$. 
\item Suppose $\G$ is intransitive. Let $\mbox{Rep}(G)$ be the monoidal 
category of representations of $G$. Then 
$\Hom(\cV_{\xi},\cV_{\psi})=\mbox{Rep}(G)$, if $\xi=\psi$, and zero otherwise. 
\item Suppose $\G$ is free, so that $E$ is a normal subgroup of $G$. 
Note that in this case $\xi$ (and $\phi$) have to be characters which are 
constant on each conjugacy class of $E$. Theorem~\ref{catchar} says that the 
restriction of $\tilde{h}$ to $E$ is the scalar representation of $G$ 
determined by $\psi\xi^{-1}$. A concrete example of some interest is the case 
where $E=\{\pm 1\}\subset SU(2)=G$. We have two characters 
on $E$, namely the trivial one, say $\xi$, and the inclusion 
$\{\pm 1\}\subset\bbC^*$, say $\psi$. Now an easy exercise reveals 
that we have 
{\em
$$\Hom(\xi,\xi)=\Hom(\psi,\psi)=\mbox{Rep}(SU(2))_{\scriptstyle{\mbox{even}}}=
\mbox{Rep}(SO(3))$$ 
}
and 
{\em
$$\Hom(\xi,\psi)=\Hom(\psi,\xi)=\mbox{Rep}(SU(2))_{\scriptstyle
{\mbox{odd}}},$$
}
where the last expression denotes the subcategory of {\em$\mbox{Rep}(SU(2))$} 
generated by the odd spins only. 
\end{enumerate}
\end{Ex}

\section{Concluding remarks}

In this section we give a rather incomplete sketch of some features 
of the categorical representations of general categorical groups and also 
make some remarks about possible generalizations of our constructions. 
 
As already proved, a (strict) categorical representation corresponds 
precisely to a homomorphism of crossed modules
\begin{equation}
\label{catrep2}
\begin{CD}
E(\G)@>{\partial}>>G(\G)\\
@V{R_p}VV @VV{R_b}V\\
(\bbC^*)^N@>{1}>> S_N.
\end{CD}
\end{equation}
This leads to the following concrete description of a categorical 
representation.

\begin{Lem} An $N$-dimensional categorical representation of $\G$ consists 
of a group homomorphism $R_b\colon G(\G)=\G_0\to S_N$ whose kernel contains 
the image of $\partial$, together with $N$ group 
characters, $\xi_1,\ldots,\xi_n$, on $E=E(\G)$ satisfying 
\begin{equation}
\label{catrep3}
\xi_i(X\triangleright e)=\xi_{R_b(X)i}(e),
\end{equation}
for any $X\in G(\G)$, $e\in E$ and $i=1,\ldots,N$. Here $R_b(X)$ denotes the left 
action on the set of $N$ elements.  
\end{Lem}
\noindent{\bf Proof} Both $R_b$ and $R_p$ in (\ref{catrep2}) are group 
homomorphisms. Clearly $R_p$ is a group homomorphism if and only if it 
defines $N$ group characters. Also (\ref{catrep2}) is 
commutative if and only if the kernel of $R_b$ contains the image of 
$\partial$. Finally, (\ref{catrep3}) is equivalent to the action 
condition on $(R_b,R_p)$. 
\qed

\begin{Rem} Note that (\ref{catrep3}) implies that each $\xi_i$ has to be 
invariant under the action of $\ker(R_b)$, i.e. 
$$\xi_i(X\triangleright e)=\xi_i(e),$$
for any $X\in \ker(R_b)$ and $e\in E$.
\end{Rem}

\begin{Def} Let $R$ be a categorical representation. 
\begin{description}
\item{a)} $R$ is called {\em decomposable} 
if there are two non-zero categorical representations $S$ and $T$ such that 
$R\cong S\boxplus T$ holds. $R$ is called {\em indecomposable} if it is not 
decomposable.
\item{b)} $R$ is called {\em reducible} if there exist a categorical 
representation $S$ of dimension less than $R$ and two 1-intertwiners $h\colon S\to R$ and 
$h'\colon R\to S$ such that $h\circ h'\colon S\to S$ is isomorphic to 
the identity 1-intertwiner on $S$. $R$ is called {\em irreducible} if it 
is not reducible. 
\end{description}
\end{Def}

Note that any irreducible categorical representation is indecomposable as 
well, but the converse is false as the following example shows. 

\begin{Ex} In the example $\G(2,3)$ of the previous section, we saw that 
$\V(2)\cong\V(1)\boxplus\V(1)$ is decomposable and $\V(2)_{\xi}$ is 
irreducible for $\xi\ne 1$. For $\xi=1$, we see that $\V(2)_1$ is 
indecomposable, but reducible. The lemma below shows that all 
categorical representations of $\G(2,3)$ of dimension greater than two are 
decomposable.
\end{Ex}

\begin{Rem} The appearance of indecomposable reducible categorical 
representations is of course due to the fact that the 1-morphisms in 2Vect 
are matrices with only non-negative integer entries. For example, to decompose 
the non-trivial 2-dimensional representation $R_b\colon C_2\to S_2$ as a 
representation of groups, one uses an intertwiner with negative and fractional 
entries. Consequently the monoidal 2-category of categorical representations 
of $\G(2,3)$ is not semi-simple (see~\cite{M-S} for the precise definition of 
semi-simplicity). Therefore it is not clear if it can be used for the 
construction of topological state-sums, because semi-simplicity is an 
essential ingredient in the proof of topological invariance of these 
state-sums. Note that one cannot simply ignore the indecomposable reducible 
categorical representations, because, as we showed in the previous section, 
one of the summands in the decomposition of 
$\V(2)_{\xi}\boxtimes\V(2)_{\psi}$ is equal to $\V(2)_1$, for any 
choice of $\xi$ and $\psi$. It seems that this problem also exists in 
Crane and Yetter's generalization of categorical representations~\cite{CY-MC}.  
\end{Rem}

There are two elementary results about indecomposable categorical 
representations that we decided to include in these remarks because they 
are very easy to prove.  

\begin{Lem} An $N$-dimensional categorical representation, $R$, is 
indecomposable if and only if the action of $R_b(G(\G))$ on 
$\{1,\ldots,N\}$ is transitive. 
\end{Lem}

\noindent{\bf Proof} Suppose we have $R\cong S\boxplus T$ and $\dim S=K$ and 
$\dim T=N-K$. Clearly we can write 
$R_b\cong S_b\oplus T_b$, where $\oplus$ means the composite of 
$S_b\times T_b$ and the canonical map between the symmetric groups 
$S_K\times S_{N-K}\to S_{K+(N-K)}=S_N$, and $\cong$ means equal up to 
conjugation by a fixed permutation. Thus we see that the action of 
$R_b(G(\G))$ on $\{1,\ldots,N\}$ is not transitive. 

Conversely, suppose that $\{1,\ldots,N\}$ can be written as the union of two 
non-empty subsets $A$ and $B$ which are both invariant under the action of 
$R_b(G(\G))$. By reordering we may assume that $A=\{1,\ldots,K\}$ and 
$B=\{K+1,\ldots,N\}$. The restrictions of $R_b(G(\G))$ to $A$ and $B$ 
respectively yield a decomposition $R_b=S_b\oplus T_b$. Take $S_p=
(\xi_1,\ldots,\xi_K)$ and $T_p=(\xi_{K+1},\ldots,\xi_N)$. Then condition 
(\ref{catrep3}) shows that $(S_b,S_p)$ and $(T_b,T_p)$ are both categorical 
representations and $R=S\oplus T$. Note that $\cong$ had become the identity 
here because of our assumption that $A=\{1,\ldots,K\}$ and 
$B=\{K+1,\ldots,N\}$, which corresponds to the choice of a fixed permutation. 
\qed

\begin{Lem}
If $R$ is indecomposable, then it is completely determined by 
$R_b$ and just one character on $E(\G)$ which is invariant 
under the action of $\ker(R_b)$. 
\end{Lem}
\noindent{\bf Proof} Suppose $R$ is indecomposable. 
Let $i\ne j\in\{1,\ldots,N\}$ be arbitrary. 
By assumption there exists an $X\in G(\G)$ such 
that $R_b(X)i=j$. By (\ref{catrep3}) we 
have $\xi_i(X\triangleright f)=\xi_{j}(f)$. Thus, if you fix $\xi_i$, 
then $\xi_j$ is uniquely determined. Since $i,j$ were 
arbitrary, this shows that one character determines uniquely all the others. 
We already remarked that any such character has to be invariant under 
$\ker(R_b)$.
\qed

\begin{Rem} Given a subgroup $H\subset G=G(\G)$ of index $N$, the action of 
$G$ on $G/H$ by left (or right) multiplication is transitive. Given an 
ordering on the elements of $G/H$ we thus obtain a group homomorphism 
$R_b\colon G\to S_N$ such that $R_b(G)$ acts transitively on 
$\{1,\ldots,N\}$. Another choice of ordering leads to a conjugate group 
homomorphism. For any $X\in G$, the construction above applied to 
$XHX^{-1}$ yields the homomorphism $XR_bX^{-1}$. It is easy to check that 
the converse is also true: given $R_b$ 
satisfying the above condition, the kernel of $R_b$ has index $N$ in 
$G$ and any group homomorphism 
$G\to S_N$, determined by an ordering on $G/\ker(R_b)$, is 
conjugate to $R_b$. This sets up a bijective 
correspondence between conjugacy classes of group homomorphisms $G\to S_N$ 
and conjugacy classes of subgroups of $G$ of index $N$.

There is quite some literature on the theory of permutation 
representations of finite groups. 
In this theory the building blocks are the transitive permutation 
representations. The key observation about them, from the point of view of representation theory, is 
that Mackey's theory of induced representations carries over to the context of 
transitive permutation representations without problems \cite{LLC,LLBC}. This allows for a complete description of the decomposition of the tensor product of two 
transitive permutation representations into a direct sum of 
transitive permutation representations, for example. Clearly categorical 
representations, as defined in this paper, are a generalization of permutation 
representations, the indecomposable ones being the generalizations of the transitive 
permutation representations, and one could probably 
generalize some of the techniques used for the latter to study 
categorical representations.       
\end{Rem}   
 
Although we have not worked out all the details about general categorical 
representations, we can say some more about the categorical representations 
of intransitive categorical groups $\G$, of which $\G(2,3)$ is a very simple 
example. Despite its simplicity the case of $\G(2,3)$ reveals the general 
features: an indecomposable categorical representation $\V$ corresponds to 
an orbit, $G(\G)\xi$, in the set of characters of $E(\G)$,  
and a homomorphism $G(\G)\to S_N$ corresponding to a subgroup 
$1\subseteq H\subseteq G(\G)_{\xi}$, where $G(\G)_{\xi}$ is the stabilizer 
of $\xi$ in $G(\G)$. $\V$ is irreducible if and only $H=G(\G)_{\xi}$. 
1-Intertwiners between two categorical representations can be interpreted as 
homogeneous vector bundles on the cartesian product of 
the two orbits in the spaces of characters of $E(\G)$. 2-Intertwiners 
can be interpreted as maps between homogeneous bundles. The decomposition of 
the monoidal product of two indecomposable categorical representions $\V_1$ 
and $\V_2$ can be obtained by looking at the decomposition into orbits of 
the cartesian product of the two orbits, corresponding to $\V_1$ and $\V_2$ 
respectively, and can, in general, contain indecomposable reducible 
categorical representations. 

A short remark on the fact that we have only worked out examples of 
strict categorical representations is in place. If $\G$ is intransitive, then 
nothing terribly interesting happens when considering general categorical 
representations. Each $N$-dimensional indecomposable categorical 
representations is determined by a homomorphism of crossed modules, just as 
in the strict case, together with an additional group 2-cocycle on $G(\G)$ 
with values in $(\bbC^*)^N$. Because we have only considered strict units, 
the 2-cocycles have to be normalized, but that is the only restriction. The 
1-intertwiners between two indecomposable categorical representations 
$\V_1$ and $\V_2$ become projective homogeneous vector bundles on the 
cartesian product of the orbits corresponding to $\V_1$ and $\V_2$, 
with a projective action of $G(\G)$ which fails to be an ordinary action by 
the quotient of the 2-cocycles of $\V_1$ and $\V_2$. The 2-intertwiners are 
just maps between these vector bundles which intertwine the projective 
actions. However, for general $\G$ the categorical representations are not 
so easy to describe. We tried to weaken the notion of homomorphism between 
crossed modules by introducing a 2-cocycle, but we found that this 2-cocycle 
has to satisfy an additional independent equation. Because the interpretation 
of this equation is not clear, we decided to work out the 
strict categorical representations only.   

Finally we want to comment on possible generalizations of our framework. One 
can consider generalizations of the notion of categorical representation and 
generalizations of the notion of categorical group. As an example of the former 
one could consider monoidal bicategories other than 2Vect. For example, one 
could consider Crane and Yetter's monoidal bicategory of measurable categories~\cite{CY-MC}. 
As they remark in their introduction, this allows for more interesting 
categorical representations of categorical Lie groups because the 
base group can be represented in more general topological symmetry groups than 
$S_N$. However, as we remarked above, indecomposable reducible categorical 
representations seem to appear in this setting as well, which puts its 
applicability for state-sums at risk. Somehow the discreteness of the non-negative 
integer entries in the 1-morphisms in 2Vect has not been solved completely in 
this new setting. 

If one only considers 1-dimensional categorical representations, then, 
of course, no indecomposable reducible categorical representations appear, e.g. 
the last example of the previous section yields a semi-simple 
monoidal 2-category and could be used for the construction of topological state-sums 
of 3- and 4-dimensional manifolds. We have not worked out what these invariants 
are, but possibly they are connected to Yetter's~\cite{Y-EX} and Roberts~\cite{R-EX} refined invariants. One idea for a generalization would be to replace $\Vect$ 
by a more interesting braided monoidal category, such as the ones appearing in the 
representation theory of quantum groups.

\end{document}